\documentclass[10pt]{article}
\usepackage{amsmath, amssymb, braket, amsthm, graphicx, titlefoot, color}
\usepackage{setspace}
\usepackage{authblk}
\usepackage[square,numbers]{natbib}

\newtheorem{theorem}{Theorem}

\newcommand*{\affaddr}[1]{#1} 
\newcommand*{\affmark}[1][*]{\textsuperscript{#1}}

\title{Centre Manifold Analysis of 3-D nonlinear system and Kinetic stability of Protein Assembly}
\author{%
	Souma Mazumdar\affmark[1]
	and Gautam Gangopadhyay\affmark[2]\\
	\affaddr{\affmark[1]Department of Theoretical Sciences}\\
	\affaddr{\affmark[2]Department of Chemical, Biological and Macro-molecular Sciences}\\
	\affaddr{S. N. Bose National Centre for Basic Sciences \\ Block - JD, Sector - III, Salt Lake City, Kolkata - 700 106}
	\thanks{\affmark[1]Email: \texttt{souma.mazumdar@bose.res.in}}
	\thanks{\affmark[2]Email: \texttt{gautam@bose.res.in}}
}
 
\date{}

\begin{document}
	\maketitle 
	
	\begin{abstract}
		Centre Manifold analysis of a $3-D$ nonlinear system with general second order nonlinearities have been worked out. The system is shown to possess two fixed points on the reduced $2-D$ centre manifold. By introducing a $2-D$ centre manifold one can show how an oscillatory dynamics may be generated in the system.	We also state and prove a theorem to find the stability of the resultant centre manifold equation apriori from the parity of the  nonlinear terms in the  original equations. For a $2-D$ nonlinear model with the example picked up from biochemistry, the protein molecules in assembly, kinetic stability analysis is provided for the chosen example and establish herewith the validity of the theorem for our chosen example.
		
	\end{abstract}

	\section{Introduction}
    Stability analyis of vector fields around an equilibrium point is a subject of interest in nonlinear dynamical systems. 
     For the nonlinear system if the eigenvalues of the linearized Jacobian has a mixture of positive, negative and zero real parts the analysis is not so straightforward and we have to go to the centre manifold theory to properly analyse the dynamics. The eigenspaces of the eigenvectors corresponding to positive, negative and zero real parts are considered. We assign a 3-D coordinate system for representing the eigenspaces of these three particular types of eigenvectors. But in case of centre manifold theory we consider only a 2-D coordinate system for representing the eigenspaces of the eigenvectors of negative and zero real parts of eigenvalues. They are termed as the stable and centre directions along the $Y$ and $X$ axis, respectively so that  it boils down only to the analysis along the centre direction to determine stability.    In case the spectrum of the linearized Jacobian contains purely imaginary eigenvalues, a slowly decaying oscillatory amplitude can appear. An investigation with application of centre manifold theory for limit cycle behaviour\cite{liu2000application}  are also made in the context of nonlinear aeroelasticity. Limit cycle     stability reversal of a two and three degrees of freedom
    airfoils has been investigated in \cite{dessi2002limit,dessi2004limit}. Limitations of centre manifold theory for limit cycle calculation has been shown in \cite{grzkedzinski2005limitation}. \par
    
    In our work we primarily focus on application of centre manifold theory for a $3-D$ system with general second order nonlinear terms. We perform the stability analysis by reducing the state equations on a $2-D$ centre manifold. The centre manifold is found to possess two fixed points one stable and another unstable. The system moves towards the stable fixed point starting from the unstable fixed point in asymptotic time. Along the course of discussion we also state and prove a theorem along with a conjecture made in the theorem.
    Then we present an example of a system which describes the kinetics of the protein molecules\cite{sanchez2010protein,manning2004structural} in an assembly. Stability analysis is done for it in the light of centre manifold theory. \par
    The paper is organised as follows. In the prelude we introduce the centre manifold theory, its various ramifications and technical details like eigenbasis and similarity transformations, block diagonal representation of the system, seek for the leading order term through centre manifold equation. All these are discussed at length to adapt the  mathematical concepts for applied dynamical systems\cite{strogatz2018nonlinear,thompson2002nonlinear} with a few degrees of freedom. Then we discuss about a general 3-D system with second order nonlinearities for  the stability analysis.  Then we state and prove a theorem for the condition of stability and in our example of the $2-D$ model through centre manifold theory we have shown the  verification of the theorem in the light of obtained results. Finally we conclude with a short and modest concluding paragraph.
    \section{Centre Manifold Theory}
    Consider a $n\times n$ nonlinear system represented by the following vector equation.
    \begin{equation}
    \begin{split}
    \dot{\vec{z}}=\vec{f}(\vec{z})
    \end{split}
    \end{equation}
    where $f(\vec{z})$ is comprised of linear and nonlinear terms of the individual variables.\\
    To analyse the sysem using centre manifold theory we have to transform the equations into the following form\cite{wiggins2003introduction}
    \begin{equation}\label{Decouple}
    \begin{split}
    &\dot{\vec{x}}=A\vec{x}+\vec{f}(\vec{x},\vec{y})\\&
    \dot{\vec{y}}=B\vec{y}+\vec{g}(\vec{x},\vec{y})
    \end{split}
    \end{equation}
    where $\vec{x}$ is a $c$ dimensional vector and $A$ is a $c\times c$ constant matrix having eigenvalues of zero real parts. Similarly $\vec{y}$ is a $s$ dimensional vector and $B$ is a $s\times s$ constant matrix having eigenvalues of negative real part. We have $c+s=n$. Here $\vec{f}(\vec{x},\vec{y})$ and $\vec{g}(\vec{x},\vec{y})$ are the nonlinear parts of the system. The advantage of such a transformation is twofold. On one hand we decouple the centre direction from the stable direction and also separate the linear part from the nonlinear part. But if the system cannot be written in this form then we have to go for a change of basis transformation. We have to transform the system in its eigenbasis. Transforming the system in its eigenbasis may lead to interprete the system in its simplest form and brings the system in the form as represented by the equation (\ref{Decouple}). For converting to eigenbasis we follow the following procedure.
    
    \subsection{Eigenbasis analysis}
    For converting to eigenbasis we have to first compute the linearized Jacobian($J$) around the equilibrium point. Then we have to compute the eigenvalues of the Jacobian. Let there be $c$ eigenvalues with zero real parts and $s$ eigenvalues with negative real parts. The eigenvectors of the eigenvalues of the zero real part span the centre eigenspace represented by the axis $X$. The eigenvectors of the eigenvalues of negative real part span the stable eigenspace represented by $Y$ axis. Let the eigenvectors of the zero real part eigenvalues be $n_{c_{1}}$ \dots $n_{c_{c}}$. Let the eigenvectors of negative real part eigenvalues be $n_{s_{1}}$ \dots $n_{s_{s}}$. Let the vector field in the eigenbasis be represented by $U$. We do the following transformation.
   \begin{equation}
   \begin{split}
   \vec{X}=P\vec{U}
   \end{split}
   \end{equation}
   where $\vec{X}$ is the original vector field and $P$ is a $n\times n$ matrix which is constructed as follows.
   \begin{equation}
   \begin{split}
   P=\begin{bmatrix}
   n_{c_{1}} \dots n_{c_{c}} \; n_{s_{1}} \dots n_{s_{s}}
   \end{bmatrix}
   \end{split}
   \end{equation} 
   If $\dot{\vec{X}}=E\vec{X}+\vec{f}(\vec{x},\vec{y})$, where $\vec{f}$ represents the nonlinear part and thereby the original system of equations on transformation of $\vec{X}=P\vec{U}$ the equation takes the form,
   \begin{equation}\label{Transformed vector eqn}
   \begin{split}
   & P\dot{\vec{U}}=EP\vec{U}+\vec{f}(\vec{x},\vec{y}) \\ &
   \implies \dot{\vec{U}}=P^{-1}EP\vec{U}+P^{-1}\vec{f}(\vec{x},\vec{y}).
   \end{split}
   \end{equation}
   The transformation of $E \rightarrow P^{-1}EP$ is known as similarity transformation. Now it can be checked the matrix $P^{-1}EP$ is in block diagonal form and in case of real eigenvalues it will be a diagonal matrix with all the eigenvalues(both zero and negative) aligned along the diagonal. The system can be reduced to the form as given by equation (\ref{Decouple}) with centre directions decoupled from the stable directions. If it is a two dimensional system with $n=2$ we should have $c=1$ and $s=1$. Then the matrix $A$ in equation (\ref{Decouple}) will be just the eigenvalue of the centre direction that is $0$ and matrix $B$ will just be the eigenvalue of the stable direction that is a negative number.
    
    \subsection{Centre Manifold analysis}
    Once the system is in the form as given by equation (\ref{Decouple}) we can carry forward our analysis from that point.\par
    Consider the decoupled directions $\vec{x}$ and $\vec{y}$. We write,
    \begin{equation}
    \begin{split}
    \vec{y}=\vec{h}(\vec{x})
    \end{split}
    \end{equation}
    where $\vec{h}$ represent the expressions we use to represent $\vec{y}$ in terms of $\vec{x}$. We call the above equation is the centre manifold equation. We consider $\vec{h}(0)=0$ and $D\vec{h}(0)=0$ where $0$ is the equilibrium point.The first condition implies that $y=0$ at the equilibrium point and the second condition implies that the centre manifold  is tangent to the centre eigenspace at the point $0$.  We expand $\vec{h}(x)$ in Taylor series expansion. Then we can write,
    \begin{equation}
    \vec{y}=\vec{h}(\vec{x})=a_{0}\vec{x^2}+a_{1}\vec{x^3}+\mathcal{O}(4)
    \end{equation}
    Then,
    \begin{equation}
    \dot{\vec{y}}=D\vec{h}(\vec{x})\dot{\vec{x}}=(2a_{0}\vec{x}+3a_{1}\vec{x^2}+\mathcal{O}(3))\dot{\vec{x}}
    \end{equation}
    where $D$ stands for the jacobian of the vector $\vec{h}$. For $2$ dimensional system $D\vec{h}$ will be just the derivative of function $h$ with respect to $x$. The above equation can be written as,
    \begin{equation}\label{N(x)=0}
    \begin{split}
    &\mathcal{N}(\vec{x})=(2a_{0}\vec{x}+3a_{1}\vec{x^2}+\mathcal{O}(3))\dot{\vec{x}}-\dot{\vec{y}}=0 \\ &
    \implies \mathcal{N}(\vec{x})=(2a_{0}\vec{x}+3a_{1}\vec{x^2}+\mathcal{O}(3))(A\vec{x}+\vec{f}(\vec{x},\vec{y}))-(B\vec{y}+g(\vec{x},\vec{y}))=0 \\ &
    \implies \mathcal{N}(\vec{x})=(2a_{0}\vec{x}+3a_{1}\vec{x^2}+\mathcal{O}(3))(A\vec{x}+\vec{f}(\vec{x},\vec{h}(\vec{x})))-(B\vec{h}(\vec{x}))+g(\vec{x},\vec{h}(\vec{x})))=0
    \end{split}
    \end{equation}
    Now $\mathcal{N}(\vec{x})$ is a system of equations containing only $\vec{x}$. For a two dimensional system $\mathcal{N}$ will be just a expression in  a single variable $x$. We can equate the coefficients of various powers of $x$ to zero and solve for $a_{0}$, $a_{1}$ etc. We will always consider the leading order term. If $a_{0}=0$ then we go to the next higher order term and solve for $a_{1}$ and so on. Then the centre manifold will be given by the equation,
     \begin{equation}
   \vec{y}=\vec{h}(\vec{x})=a_{0}\vec{x^2}+a_{1}\vec{x^3}+\mathcal{O}(4)
   \end{equation}
   upto leading order. Once we obtain $\vec{y}$ in terms of $\vec{x}$ that is the centre manifold equation we no longer consider the $\vec{y}$ direction. We consider only the $\vec{x}$ direction and write,
   \begin{equation}
   \dot{\vec{x}}=A\vec{x}+\vec{f}(\vec{x},\vec{h}(\vec{x}))
   \end{equation}
   The above equation is called the reduced system of equations. We can use it to analyse stability around the equilibrium point.
   
   \subsection{Centre Manifold Analysis of a general 3-D system of second order nonlinearities}
   
  Let the linearized Jacobian around the equilibrium point yields eigenvalues $i\lambda_{1},-i\lambda_{1},-\lambda_{2}$.
  We consider $\lambda_{1} > 0$ and $\lambda_{2} > 0$. The real parts of first two eigenvalues is $0$, so there is a centre manifold. The third eigenvalue is real negative, so there is a stable manifold. Let on eigenbasis transformation the system equations take the form
  \begin{equation}
  \begin{split}
  \begin{bmatrix}
  \dot{x} \\ \dot{y} \\ \dot{z}
  \end{bmatrix}=\begin{bmatrix}
  0 & \lambda_{1} & 0 \\
  -\lambda_{1} & 0 & 0 \\
  0 & 0 & -\lambda_{2}
  \end{bmatrix}\begin{bmatrix}
  x \\ y \\z 
  \end{bmatrix}+\begin{bmatrix}
  f_{1}(x,y,z) \\ f_{2}(x,y,z) \\ f_{3}(x,y,z)
  \end{bmatrix}
  \end{split}
  \end{equation}
  where $f_{1},f_{2},f_{3}$ are the nonlinear parts containing second order nonlinearities in $x,y,z$.  We can assume some second order polynomial structure in $x,y,z$ for $f_{1},f_{2},f_{3}$. Let,
  \begin{equation}
  \begin{split}
  & f_{1}(x,y,z)=c_{0}x^2+c_{1}y^2+c_{2}z^2+c_{3}xy+c_{4}yz+c_{5}zx \\ &
  f_{2}(x,y,z)=d_{0}x^2+d_{1}y^2+d_{2}z^2+d_{3}xy+d_{4}yz+d_{5}zx \\ &
  f_{3}(x,y,z)=e_{0}x^2+e_{1}y^2+e_{2}z^2+e_{3}xy+e_{4}yz+e_{5}zx
  \end{split}
  \end{equation}
  Therefore the system equations are given by,
  \begin{equation}
  \begin{split}
  & \dot{x}=\lambda_{1}y+c_{0}x^2+c_{1}y^2+c_{2}z^2+c_{3}xy+c_{4}yz+c_{5}zx \\ &
  \dot{y}=-\lambda_{1}x+d_{0}x^2+d_{1}y^2+d_{2}z^2+d_{3}xy+d_{4}yz+d_{5}zx \\ &
  \dot{z}= -\lambda_{2}z+e_{0}x^2+e_{1}y^2+e_{2}z^2+e_{3}xy+e_{4}yz+e_{5}zx
  \end{split}
  \end{equation}
  The equilibrium point is the origin$(0,0,0)$. If the equilibrium point is a point other than the origin say $(x_{0},y_{0},z_{0})$ then we would have made the transformation $x \rightarrow x-x_{0}, y \rightarrow y-y_{0}, z \rightarrow z-z_{0}$ and have brought the equilibrium position to the origin. \\ The centre manifold is two dimensional and the stable manifold is one dimensional. Let the centre manifold is given by the following taylor series expansion.
  \begin{equation}
  z=h(x,y)=a_{0}x^2+a_{1}xy+a_{2}y^2+\mathcal{O}(3)
  \end{equation}
  We take upto $\mathcal{O}(2)$. So
  \begin{equation}
  z=h(x,y)=a_{0}x^2+a_{1}xy+a_{2}y^2
  \end{equation}
  Therefore,
  \begin{equation}
  \begin{split}
  & \dot{z} =\frac{\partial h}{\partial x}\dot{x}+\frac{\partial h}{\partial y}\dot{y} \\ &
  \implies \dot{z}-(2a_{0}x+a_{1}y)\dot{x}-(2a_{2}y+a_{1}x)\dot{y}=0
  \end{split}
  \end{equation}
  Substituting the values of $\dot{x},\dot{y},\dot{z}$ from the system equations and replacing $z$ by $(a_{0}x^2+a_{1}xy+a_{2}y^2)$ in the equations we have,
  \begin{equation}
  \begin{split}
  & -\lambda_{2}(a_{0}x^2+a_{1}xy+a_{2}y^2)+e_{0}x^2+e_{1}y^2+e_{2}(a_{0}x^2+a_{1}xy+a_{2}y^2)^2+e_{3}xy\\ & +e_{4}y(a_{0}x^2+a_{1}xy+a_{2}y^2)+e_{5}x(a_{0}x^2+a_{1}xy+a_{2}y^2)\\ &-(2a_{0}x+a_{1}y)[\lambda_{1}y+c_{0}x^2+c_{1}y^2+c_{2}(a_{0}x^2+a_{1}xy+a_{2}y^2)^2+c_{3}xy \\ &+c_{4}y(a_{0}x^2+a_{1}xy+a_{2}y^2)+c_{5}x(a_{0}x^2+a_{1}xy+a_{2}y^2)]\\ &-(2a_{2}y+a_{1}x)[-\lambda_{1}x+d_{0}x^2+d_{1}y^2+d_{2}(a_{0}x^2+a_{1}xy+a_{2}y^2)^2+d_{3}xy \\ & +d_{4}y(a_{0}x^2+a_{1}xy+a_{2}y^2)+d_{5}x(a_{0}x^2+a_{1}xy+a_{2}y^2)]=0
  \end{split}
  \end{equation}
  
  Considering only $\mathcal{O}(2)$ and equating the coefficients of $x^2,xy,y^2$ to $0$ respectively we have,
  \begin{equation}
  \begin{split}
  &-\lambda_{2}a_{0}+e_{0}+a_{1}\lambda_{1}=0 \\ &
  -\lambda_{2}a_{1}+e_{3}-2a_{0}\lambda_{1}+2a_{2}\lambda_{1}=0 \\ &
  -\lambda_{2}a_{2}+e_{1}-a_{1}\lambda_{1}=0
  \end{split}
  \end{equation}
  Solving for $a_{0},a_{1},a_{2}$ we have,
  \begin{equation}
  \begin{split}
 & a_{0}=\frac{\lambda_{2}(e_{0}\lambda_{2}+e_{3}\lambda_{1})+4e_{0}\frac{\lambda_{1}^2}{\lambda_{2}}+2\lambda_{1}^2(e_{1}-e_{0})}{\lambda_{2}^3+4\lambda_{1}^2} \\ &
 a_{1}=\frac{\lambda_{2}^2e_{3}+2\lambda_{1}\lambda_{2}(e_{1}-e_{0})}{\lambda_{2}^3+4\lambda_{1}^2} \\ &
 a_{2}=\frac{\lambda_{2}(e_{1}\lambda_{2}-e_{3}\lambda_{1})+4e_{1}\frac{\lambda_{1}^2}{\lambda_{2}}+2\lambda_{1}^2(e_{0}-e_{1})}{\lambda_{2}^3+4\lambda_{1}^2}
  \end{split}
  \end{equation}
  Therefore the centre manifold equation is given by,
  \begin{equation}
  \begin{split}
 & z=\frac{\lambda_{2}(e_{0}\lambda_{2}+e_{3}\lambda_{1})+4e_{0}\frac{\lambda_{1}^2}{\lambda_{2}}+2\lambda_{1}^2(e_{1}-e_{0})}{\lambda_{2}^3+4\lambda_{1}^2}x^2+\frac{\lambda_{2}^2e_{3}+2\lambda_{1}\lambda_{2}(e_{1}-e_{0})}{\lambda_{2}^3+4\lambda_{1}^2}xy
 \\ & +\frac{\lambda_{2}(e_{1}\lambda_{2}-e_{3}\lambda_{1})+4e_{1}\frac{\lambda_{1}^2}{\lambda_{2}}+2\lambda_{1}^2(e_{0}-e_{1})}{\lambda_{2}^3+4\lambda_{1}^2}y^2
  \end{split}
  \end{equation}
  The reduced equations on the centre manifold on considering terms of $\mathcal{O}(2)$ and neglecting higher order terms is given by,
  \begin{equation}\label{reduced eqn x y}
  \begin{split}
 & \dot{x}=\lambda_{1}y+c_{0}x^2+c_{1}y^2+c_{3}xy \\ &
 \dot{y}=-\lambda_{1}x+d_{0}x^2+d_{1}y^2+d_{3}xy
  \end{split}
  \end{equation}
  Now we seek to find if there is any oscillations on this reduced phase space.
  For this we write,
  \begin{equation}\label{r theta coord}
  \begin{split}
  &x=r\cos\theta \\ &
  y=r\sin\theta
  \end{split}
  \end{equation}
  on transformation to $r,\theta$ coordinate. Then
  \begin{equation}
  r^2=x^2+y^2
  \end{equation}
  Differentiating both sides with respect to time $t$ we have,
  \begin{equation}
  \begin{split}
 & r\dot{r}=x\dot{x}+y\dot{y} \\ &
 \implies \dot{r}=\frac{x\dot{x}+y\dot{y}}{r}
  \end{split}
  \end{equation}
  On puting the values of $\dot{x}$, $\dot{y}$ and simplifying we have,
  \begin{equation}
  \begin{split}
  \dot{r}=\frac{x(c_{0}x^2+c_{1}y^2+c_{3}xy)+y(d_{0}x^2+d_{1}y^2+d_{3}xy)}{r}
  \end{split}
  \end{equation}
  At this point for the sake of simplification we assume values for $c_{0},c_{1},c_{3},d_{0},d_{1},d_{3}$.\\  We consider $c_{0}=c_{1}=c_{3}=d_{0}=d_{1}=d_{3}=1$. \\ Then the above equation takes the form,
  \begin{equation}\label{r dot}
  \begin{split}
  \dot{r}& =\frac{x(x^2+y^2+xy)+y(x^2+y^2+xy)}{r} \\  &
  \implies \dot{r}=r^2(\cos\theta+\sin\theta)(1+\cos\theta\sin\theta)
  \end{split}
  \end{equation}
  on puting $x=r\cos\theta,y=r\sin\theta$.
  Now for fixed point prediction we should solve for $r$ from $\dot{r}=0$. Therefore we have to solve for $r$ from the equation,
  \begin{equation}
  r^2(\cos\theta+\sin\theta)(1+\cos\theta\sin\theta)=0
  \end{equation}
  Now equating each of the factors to $0$ we have,
  \begin{equation}
  \begin{split}
  &r=0 \\ & \cos\theta=-\sin\theta \\ &
  \cos\theta\sin\theta=-1
  \end{split}
  \end{equation}
  The solution $r=0$ implies a fixed point at the origin which is a sink as we will see. The third solution is discarded because $\cos\theta\sin\theta > -1$ for any $\theta$. From the second solution we have,
  \begin{equation}
  \begin{split}
 & \cos\theta=-\sin\theta \\ &
 \implies \theta=n\pi+\frac{3\pi}{4}
  \end{split}
  \end{equation}
  Puting in equation (\ref{r theta coord}),
  \begin{equation}
  \begin{split}
  & x=-\frac{r}{\sqrt{2}} \\ &
  y=\frac{r}{\sqrt{2}}
  \end{split}
  \end{equation}
  for $\theta=2n\pi+\frac{3\pi}{4}$ and
   \begin{equation}
  \begin{split}
  & x=\frac{r}{\sqrt{2}} \\ &
  y=-\frac{r}{\sqrt{2}}
  \end{split}
  \end{equation}
  for $\theta=(2n+1)\pi+\frac{3\pi}{4}$. \\
  Puting the the first solution of $x,y$ in the first equation of (\ref{reduced eqn x y}) and recalling we have considered $c_{0}=c_{1}=c_{3}=1$ we have,
  \begin{equation}
  \begin{split}
  &-\frac{\dot{r}}{\sqrt{2}}=\lambda_{1}\frac{r}{\sqrt{2}}+\frac{r^2}{2}+\frac{r^2}{2}-\frac{r^2}{2} \\ &
  \implies -\dot{r}=r(\lambda_{1}+\frac{r}{\sqrt{2}})
  \end{split}
  \end{equation}
  Now the above equation is a equation only in $r$. For solving the amplitude $r$ we equate $\dot{r}=0$. Thus we have,
  \begin{equation}
  \begin{split}
  &r(\lambda_{1}+\frac{r}{\sqrt{2}})=0 \\ & 
  \implies r=0,-\sqrt{2}\lambda_{1}
  \end{split}
  \end{equation}
  Now $r=0$ implies a sink at the origin since $\dot{r} < 0$ for $r > 0$ in the neighbourhood of origin. $r=-\sqrt{2}\lambda_{1}$ is discarded because amplitude cannot be negative. \\
  Puting the second solution of $x,y$ in the first equation of (\ref{reduced eqn x y}) we have,
  \begin{equation}\label{r dot 2nd}
  \begin{split}
  \frac{\dot{r}}{\sqrt{2}}& =-\lambda_{1}\frac{r}{\sqrt{2}}+\frac{r^2}{2}+\frac{r^2}{2}-\frac{r^2}{2} \\ &
 \implies \dot{r} =r(\frac{r}{\sqrt{2}}-\lambda_{1})
  \end{split}
  \end{equation}
  For solving for amplitude we equate $\dot{r}=0$. Thus we have,
  \begin{equation}
  \begin{split}
  & r(\frac{r}{\sqrt{2}}-\lambda_{1})=0 \\
   &\implies r=0,\sqrt{2}\lambda_{1}
  \end{split}
  \end{equation}
  $r=0$ gives a sink at the origin since $\dot{r} < 0$ for $r > 0$ in the neighbourhood of origin.  $r=\sqrt{2}\lambda_{1}$ gives a source because $\dot{r}<0$ for $r<\sqrt{2}\lambda_{1}$ and $\dot{r}>0$ for $r>\sqrt{2}\lambda_{1}$ in the neighbourhood of $r=\sqrt{2}\lambda_{1}$. We want to show as $t \rightarrow \infty$ the system approaches $r=0$ starting from $r=\sqrt{2}\lambda_{1}$.\\
  For this we make the following coordinate transformation. We write $r^{\prime}=r-\sqrt{2}\lambda_{1}$. Then $r=r^{\prime}+\sqrt{2}\lambda_{1}$. We have to do this transformation because $r=\sqrt{2}\lambda_{1}$ is the starting point of the system. Puting the value of $r$ in terms of $r^{\prime}$ in equation(\ref{r dot 2nd}) we have,
  \begin{equation}
  \begin{split}
  & \dot{r^{\prime}}=(r^{\prime}+\sqrt{2}\lambda_{1})\frac{r^{\prime}}{\sqrt{2}} \\ &
  \implies \sqrt{2}\frac{dr^{\prime}}{dt}=(r^{\prime}+\sqrt{2}\lambda_{1})r^{\prime} \\ &
  \implies \frac{\sqrt{2}\lambda_{1}}{(r^{\prime}+\sqrt{2}\lambda_{1})r^{\prime}}dr^{\prime}=\lambda_{1}dt \\ &
  \implies \frac{(r^{\prime}+\sqrt{2}\lambda_{1})-r^{\prime}}{(r^{\prime}+\sqrt{2}\lambda_{1})r^{\prime}}dr^{\prime}=\lambda_{1}dt \\ &
  \implies \left[\frac{1}{r^{\prime}}-\frac{1}{r^{\prime}+\sqrt{2}\lambda_{1}}\right]dr^{\prime}=\lambda_{1}dt
  \end{split}
  \end{equation}
  Integrating both sides we have,
  \begin{equation}
  \begin{split}
  \ln\frac{r^{\prime}k}{r^{\prime}+\sqrt{2}\lambda_{1}}=\lambda_{1}t
  \end{split}
  \end{equation}
  where $k$ is the constant of integration. Then from the above equation,
  \begin{equation}
  \begin{split}
  & \frac{r^{\prime}+\sqrt{2}\lambda_{1}}{r^{\prime}k}=e^{-\lambda_{1}t} \\ &
  \implies r^{\prime}=\frac{\sqrt{2}\lambda_{1}}{ke^{-\lambda_{1}t}-1}
  \end{split}
  \end{equation}
  Puting $t=\infty$ in the above equation we get $r^{\prime}=-\sqrt{2}\lambda_{1}$. Since $r^{\prime}=r-\sqrt{2}\lambda_{1}$ when $r^{\prime}=-\sqrt{2}\lambda_{1}$ we have $r=0$. So the system approaches $r=0$ as $t \rightarrow \infty$.
     \\
  For the $\theta$ dynamics we write,
  \begin{equation}
  \begin{split}
  \theta=\tan^{-1}\frac{y}{x}
  \end{split}
  \end{equation}
  Differentiating the above equation with respect to time $t$ we have,
  \begin{equation}
  \begin{split}
  &\dot{\theta}=\frac{x\dot{y}-y\dot{x}}{x^2+y^2} \\ &
  \implies \dot{\theta}=\frac{x(-\lambda_{1}x+x^2+y^2+xy)-y(\lambda_{1}y+x^2+y^2+xy)}{x^2+y^2}
  \end{split}
  \end{equation}
  on puting the values of $\dot{x},\dot{y}$ and recalling $c_{0}=c_{1}=c_{3}=d_{0}=d_{1}=d_{3}=1$.
  On simplication it gives,
  \begin{equation}
  \begin{split}
  \dot{\theta}& =-\lambda_{1}+x-y+\frac{xy(x-y)}{x^2+y^2} \\ &
  =-\lambda_{1}+r(\cos\theta-\sin\theta)(1+\cos\theta\sin\theta)
  \end{split}
  \end{equation}
  on puting $x=r\cos\theta,y=r\sin\theta$. Since $r$ and $\theta$ are functions of $t$, $\dot{\theta}$ is also a function of $t$.  $\dot{\theta}$ is not explicitly solvable in terms of $t$ since $r$ and $\theta$ are not explicitly solvable in terms of $t$. Still we can do some analysis and prediction. If $\dot{\theta}(0) \neq 0$ due to the nonzero angular velocity the system will exhibit a periodic spiral motion with approach towards the point $r=0$. For oscillation prediction we look into the equations.
  \begin{equation}
  \begin{split}
  &x=r\cos\theta \\ &
  y=r\sin\theta
  \end{split}
  \end{equation}
  Since $r$ is approaching towards $0$ we can say the amplitude is decreasing so in the $x-t$ and $y-t$ plane the amplitude becomes small and small as $t \rightarrow \infty$ and eventually becomes $0$.
   \\ \\
   Now we proceed to prove a theorem.
   \begin{theorem}
   	Consider the system.
   	\begin{equation}
   	\begin{split}
   	&\dot{\vec{x}}=A\vec{x}+\vec{f}(\vec{x},\vec{y})\\&
   	\dot{\vec{y}}=B\vec{y}+\vec{g}(\vec{x},\vec{y})
   	\end{split}
   	\end{equation} 
   	If the nonlinear parts $\vec{f}(\vec{x},\vec{y})$ and $\vec{g}(\vec{x},\vec{y})$ are even that is $\vec{f}(\vec{x},\vec{y})=\vec{f}(-\vec{x},-\vec{y})$ and $\vec{g}(\vec{x},\vec{y})=\vec{g}(-\vec{x},-\vec{y})$ then $\vec{h}(\vec{x})$ is even that is $\vec{h}(\vec{x})=\vec{h}(-\vec{x})$
   	
   	\begin{proof}
   		We have
   		\begin{equation}
   		\begin{split}
   		&\dot{\vec{y}}=D\vec{h}(\vec{x})\dot{\vec{x}} \\ &
   	\implies	B\vec{y}+\vec{g}(\vec{x},\vec{y})=D\vec{h}(\vec{x})(A\vec{x}+\vec{f}(\vec{x},\vec{y})) \\ &
   	\implies B\vec{h}(x)+\vec{g}(\vec{x},\vec{y})=D\vec{h}(\vec{x})A\vec{x}+D\vec{h}(\vec{x})\vec{f}(\vec{x},\vec{y})
   		\end{split}
   		\end{equation}
   
   Now, $\vec{h}(\vec{x})=a_{0}\vec{x^2}+a_{1}\vec{x^3}+\mathcal{O}(4)$. Upto leading order $\vec{h}(\vec{x})$ will be $\vec{h}(\vec{x})=a_{0}\vec{x^2}$. Then $D\vec{h}(\vec{x})=2a_{0}\vec{x}$. Then we should have
   \begin{equation}\label{thm1}
   \begin{split}
   B\vec{h}(x)+\vec{g}(\vec{x},\vec{y})=(2a_{0}\vec{x})A\vec{x}+(2a_{0}\vec{x})\vec{f}(\vec{x},\vec{y})
   \end{split}
   \end{equation}
   Puting $-\vec{x}$ at both sides of the equation we have,
   \begin{equation}
   \begin{split}
   B\vec{h}(-x)+\vec{g}(-\vec{x},-\vec{y})=(-2a_{0}\vec{x})A(-\vec{x})+(-2a_{0}\vec{x})\vec{f}(-\vec{x},-\vec{y})
   \end{split}
   \end{equation}
   Since $\vec{f}(\vec{x},\vec{y})$ and $\vec{g}(\vec{x},\vec{y})$ are even we should have,
   \begin{equation}\label{thm2}
   \begin{split}
   B\vec{h}(-x)+\vec{g}(\vec{x},\vec{y})=(2a_{0}\vec{x})A\vec{x}+(-2a_{0}\vec{x})\vec{f}(\vec{x},\vec{y})
   \end{split}
   \end{equation}
   Now $\vec{f}(\vec{x},\vec{y})$ being nonlinear is at least of $\mathcal{O}(2)$. Therefore $(2a_{0}\vec{x})\vec{f}(\vec{x},\vec{y})$ is at least of $\mathcal{O}(3)$ which we do not want to consider. Therefore the last terms in equations (\ref{thm1}) and (\ref{thm2}) go to $0$. Therefore from  equations (\ref{thm1}) and (\ref{thm2}) we have.
   \begin{equation}
   \begin{split}
    & B\vec{h}(x)+\vec{g}(\vec{x},\vec{y})=(2a_{0}\vec{x})A\vec{x} \\ &
     B\vec{h}(-x)+\vec{g}(\vec{x},\vec{y})=(2a_{0}\vec{x})A\vec{x}
   \end{split}
   \end{equation}
   Subtracting one from the other we have
   \begin{equation}
   \begin{split}
  & B(\vec{h}(\vec{x})-\vec{h}(-\vec{x}))=0 \\ &
  \implies \vec{h}(\vec{x})=\vec{h}(-\vec{x})
   \end{split}
   \end{equation}
   If in $\vec{h}(\vec{x})$, $a_{0}=0$ we go to the next leading order term and consider $\vec{h}(\vec{x})=a_{1}\vec{x^3}$. Then we have $D\vec{h}(\vec{x})=3a_{1}\vec{x^2}$. Therefore in the equation $B\vec{h}(\vec{x})+\vec{g}(\vec{x},\vec{y})=D\vec{h}(\vec{x})A\vec{x}+D\vec{h}(\vec{x})\vec{f}(\vec{x},\vec{y})$ we put $\vec{h}(\vec{x})=a_{1}\vec{x^3}$. Therefore we have,
   \begin{equation}\label{thm3}
   \begin{split}
   B(a_{1}\vec{x^3})+\vec{g}(\vec{x},\vec{y})=(3a_{1}\vec{x^2})A\vec{x}+(3a_{1}\vec{x^2})\vec{f}(\vec{x},\vec{y})
   \end{split}
   \end{equation}
   Puting $\vec{x}=-\vec{x}$ in the above equation we have,
   \begin{equation}\label{thm4}
   \begin{split}
  & Ba_{1}(-\vec{x^3})+\vec{g}(-\vec{x},-\vec{y})=(3a_{1}\vec{x^2})A(-\vec{x})+(3a_{1}\vec{x^2})\vec{f}(-\vec{x},-\vec{y}) \\ &
  \implies -B(a_{1}\vec{x^3})+\vec{g}(\vec{x},\vec{y})=-(3a_{1}\vec{x^2})A\vec{x}+(3a_{1}\vec{x^2})\vec{f}(\vec{x},\vec{y})
   \end{split}
   \end{equation}
   Now $\vec{f}(\vec{x},\vec{y})$ being nonlinear is at least of $\mathcal{O}(2)$. So $(3a_{1}\vec{x^2})\vec{f}(\vec{x},\vec{y})$ is at least of $\mathcal{O}(4)$ which we do not want to consider in equations (\ref{thm3}) and (\ref{thm4}). Therefore from equations (\ref{thm3}) and (\ref{thm4}) we have,
   \begin{equation}
   \begin{split}
   & B(a_{1}\vec{x^3})+\vec{g}(\vec{x},\vec{y})=(3a_{1}\vec{x^2})A\vec{x} \\ &
   -B(a_{1}\vec{x^3})+\vec{g}(\vec{x},\vec{y})=-(3a_{1}\vec{x^2})A\vec{x}
   \end{split}
   \end{equation}
   Adding these two equations we have,
   \begin{equation}
   \begin{split}
   & 2\vec{g}(\vec{x},\vec{y})=0 \\ &
   \implies \vec{g}(\vec{x},\vec{y})=0
   \end{split}
   \end{equation}
   which is absurd since $\vec{g}(\vec{x},\vec{y})$ is a defined function in the problem. So we can say the leading order term in $\vec{h}(\vec{x})$ cannot be $\vec{x^3}$ if coefficient of $\vec{x^2}$ is $0$. So we go to the next leading order term that is $\vec{x^4}$. Considering $\vec{h}(\vec{x})=a_{2}\vec{x^4}$ and following the procedure as shown for the case when the leading order term is of $\mathcal{O}(2)$, we can show $\vec{h}(\vec{x})=\vec{h}(-\vec{x})$. We can make the conjecture that the leading order term in $\vec{h}(\vec{x})$ will be always of even power of $\vec{x}$ if the nonlinear parts are even functions of $\vec{x}$ and $\vec{y}$.
   	\end{proof}
   \end{theorem}
  In the example that we consider to illustrate the centre manifold problem we show that our above theorem and conjecture are true.
  
  \section{Kinetic Stability analysis of protein assembly around a fixed point}
  
  We pick up the nonlinear equations as given in the article\cite{tsuruyama2017kinetic}. The Nonlinear equations describe the kinetics of protein molecules. We write down the equation as follows.
  \begin{equation}\label{master eqn}
  \begin{split}
 & \dot{x}=-(D_{1}-aX_{e}-k)x+(-bX_{e}+p)z+ax^2-bxz \\ &
  \dot{z}=(k-cY_{e})x-pz+cx^2+cxz
  \end{split}
  \end{equation}
  The above two equations form the master equations for the analysis that follows. Following the path as shown in the paper\cite{tsuruyama2017kinetic} we consider sufficiently small values of $D_{1},k,Y_{e},p$. They can be set to $0$. Then the equations (\ref{master eqn}) can be written down as follows.
  \begin{equation}
  \begin{split}
 & \dot{x}=aX_{e}x-bX_{e}z+ax^2-bxz \\ &
 \dot{z}=cx^2+cxz
  \end{split}
  \end{equation}
  The equilibrium point is $(0,0)$. The linearized Jacobian$(J)$ around the equilibrium point $(0,0)$ is givven by
  \begin{equation}
  J=\begin{bmatrix}
  aX_{e} & -bX_{e} \\  0 & 0
  \end{bmatrix}
  \end{equation}
  The eigenvalues of $J$ are $\lambda=0,aX_{e}$. So we see one of the eigenvalues is 0 which indicates that it has a centre eigenspace and permits us to do the centre manifold anaysis. The other eigenspace is a stable eigenspace corresponding to negative eigenvalue. This prompts us to consider $aX_{e} < 0$ which implies $a<0$ because $X_{e}$ cannot be negative. The eigenvector correponding to $0$ eigenvalue is $\begin{bmatrix}
  b \\ a \end{bmatrix}$ and the eigenvector corresponding to $aX_{e}$ eigenvalue is $\begin{bmatrix} 1 \\ 0 \end{bmatrix}$. As the Jacobian $J$ is not block diagonal we have to resort to eigenbasis transformation. As mentioned in the theory we consider a transformation $\vec{X}=P\vec{U}$ where $\vec{X}$ is the original vector and $\vec{U}$ is the vector in the changed basis. In our case the system being two dimensional $\vec{X}$ and $\vec{U}$ are two dimensional vectors and $P$ is a $2 \times 2$ matrix formed by the eigenvectors of the Jacobian. We place the eigenvector of the centre direction as the first vector and that of the stable direction as the second vector. So,
  \begin{equation}
  \begin{split}
  P=\begin{bmatrix}
  b & 1 \\ a & 0
  \end{bmatrix}
  \end{split}
  \end{equation}
  The inverse is given by
  \begin{equation}
  \begin{split}
  P^{-1}=\begin{bmatrix}
  0 & \frac{1}{a} \\ 1 & -\frac{b}{a}
  \end{bmatrix}
  \end{split}
  \end{equation}
  Following equation (\ref{Transformed vector eqn}) we have
  \begin{equation}
  \dot{\vec{U}}=P^{-1}AP\vec{U}+P^{-1}\vec{f}(x,y)
  \end{equation}
  Now as mentioned in the theory $P^{-1}AP$ will be a block diagonal matrix which in our case will be just a diagonal matrix with the eigenvalues along the diagonal. We make the first diagonal entry as the eigenvalue of the centre direction and the second diagonal entry as the eigenvalue of the stable direction. Therefore we have,
  \begin{equation}
  P^{-1}AP=\begin{bmatrix} 0 & 0 \\ 0 & aX_{e} \end{bmatrix}
  \end{equation}
  Therefore,
  \begin{equation}\label{transformed basis eqn}
  \begin{split}
  \dot{\vec{U}}=\begin{bmatrix} \dot{u} \\ \dot{v} \end{bmatrix}=\begin{bmatrix} 0 & 0 \\ 0 & aX_{e} \end{bmatrix}\begin{bmatrix} u \\ v \end{bmatrix}+\begin{bmatrix} 0 & \frac{1}{a} \\ 1 & -\frac{b}{a} \end{bmatrix}\begin{bmatrix} ax^2-bxz \\ cx^2+cxz \end{bmatrix} 
  \end{split}
  \end{equation}
  Now the nonlinear part in the above equation is still in terms of $x,z$ which we would like to convert in terms of $u,v$ using the equation $\vec{X}=P\vec{U}$. Therefore we have
  \begin{equation}
  \begin{split}
  \begin{bmatrix} x \\ z \end{bmatrix}=\begin{bmatrix} b & 1 \\ a & 0 \end{bmatrix}\begin{bmatrix} u \\ v \end{bmatrix}
  \end{split}
  \end{equation}
  Therefore we have,
  \begin{equation}
  \begin{split}
 & x=bu+v \\ &
 z=au
  \end{split}
  \end{equation}
  Puting in equation (\ref{transformed basis eqn}) we have
  \begin{equation}
  \begin{bmatrix} \dot{u} \\ \dot{v} \end{bmatrix}=\begin{bmatrix} 0 & 0 \\ 0 & aX_{e} \end{bmatrix}\begin{bmatrix} u \\ v \end{bmatrix}+\begin{bmatrix} 0 & \frac{1}{a} \\ 1 & -\frac{b}{a} \end{bmatrix}\begin{bmatrix} a(bu+v)^2-b(bu+v)au \\ c(bu+v)^2+c(bu+v)au \end{bmatrix} 
  \end{equation}
  Therefore we have the system of equations as,
  \begin{equation}\label{eqns in u v}
  \begin{split}
  & \dot{u}=\frac{1}{a}\left[c(bu+v)^2+c(bu+v)au\right] \\ &
  \dot{v}=aX_{e}v+a(bu+v)^2-b(bu+v)au-\frac{b}{a}\left[c(bu+v)^2+c(bu+v)au\right]
  \end{split}
  \end{equation}
  From equation (\ref{N(x)=0}) we have
  \begin{equation}
  (2a_{0}x+3a_{1}x^2+\mathcal{O}(3))(Ax+f(x,h(x)))-(Bh(x)+g(x,h(x)))=0
  \end{equation}
  But here we replace $x$ by $u$ because we are working in the eigenbasis.\\
  Therefore we have,
  \begin{equation}\label{pre centre manifold}
  (2a_{0}u+3a_{1}u^2+\mathcal{O}(3))(Au+f(u,h(u)))-(Bh(u)+g(u,h(u)))=0
  \end{equation}
  Now we recall,
  \begin{equation}
  h(u)=a_{0}u^2+a_{1}u^3+\mathcal{O}(4)
   \end{equation}
   We consider upto $\mathcal{O}(2)$.
  Therefore,
  \begin{equation}
  \begin{split}
 & Au+f(u,h(u))=\frac{1}{a}\left[c(bu+a_{0}u^2)+c(bu+a_{0}u^2)au\right] \\ &
 Bh(u)+g(u,h(u))=aX_{e}(a_{0}u^2)+a(bu+a_{0}u^2)^2-b(bu+a_{0}u^2)au-\frac{b}{a}\left[c(bu+a_{0}u^2)^2+c(bu+a_{0}u^2)au\right]
  \end{split}
  \end{equation}
  As we consider upto $\mathcal{O}(2)$ we search for the coefficient of $u^2$ that is $a_{0}$. If $a_{0}=0$, then we go to the next higher order term that is $u^3$. Puting in equation (\ref{pre centre manifold}) we have,
  \begin{equation}
  \begin{split}
 & (2a_{0}u+3a_{1}u^2)\frac{1}{a}\left[c(bu+a_{0}u^2)+c(bu+a_{0}u^2)au\right] \\ &-\left[aX_{e}(a_{0}u^2)+a(bu+a_{0}u^2)^2-b(bu+a_{0}u^2)au-\frac{b}{a}\left[c(bu+a_{0}u^2)^2+c(bu+a_{0}u^2)au\right]\right]=0
  \end{split}
  \end{equation}
  The coefficient of $u^2$ in the above expression is $aX_{e}a_{0}-\frac{bc}{a}b^2-b^2c$. \\
  Equating the coefficient of $u^2$ to $0$ we have,
  \begin{equation}
  \begin{split}
  & aX_{e}a_{0}-b^2c(1+\frac{b}{a})=0 \\ &
  \implies a_{0}=\frac{b^2c}{aX_{e}}(1+\frac{b}{a})
  \end{split}
  \end{equation}
  \subsection{Equation of Centre manifold and reduced system}
  Equation of centre manifold is given by
  \begin{equation}
  \begin{split}
  v=h(u)=a_{0}u^2+\mathcal{O}(3)
  \end{split}
  \end{equation}
  Considering upto $\mathcal{O}(2)$ which is the leading order and substituting the value of $a_{0}$ we have,
  \begin{equation}
  v=\frac{b^2c}{aX_{e}}(1+\frac{b}{a})u^2
  \end{equation}
  The above equations is the equation of \textit{Centre Manifold}. From the first equation of the set (\ref{eqns in u v}) we have
  \begin{equation}
  \dot{u}=\frac{1}{a}\left[c(bu+v)^2+c(bu+v)au\right]
  \end{equation}
  Puting the value of $v$ as obtained in the centre manifold equation, in the above equation we have,
   \begin{equation}
  \dot{u}=\frac{1}{a}\left[c\left[bu+\frac{b^2c}{aX_{e}}(1+\frac{b}{a})u^2\right]^2+c\left[bu+\frac{b^2c}{aX_{e}}(1+\frac{b}{a})u^2\right]au\right]
  \end{equation}
  Considering upto $\mathcal{O}(2)$ of the power of $u$ in the above equation we have,
  \begin{equation}\label{reduced system}
  \dot{u}=cb(1+\frac{b}{a})u^2
  \end{equation}
  This gives the reduced system of equations.
  \subsubsection{Stability Analysis of Equilibrium Point}
  After we obtain the reduced system we no longer consider the equation of $\dot{v}$. We know the stable direction that is the $v$ direction will converge exponentially fast towards the equilibrium. So for analysis of stability of equilibrium point we can only consider equation (\ref{reduced system}).
  We have,
  \begin{equation}
  \dot{u}=cb(1+\frac{b}{a})u^2
  \end{equation}
  Now consider the neigbourhood of equilibrium point $(0,0)$ along the $u$ direction. Two cases may arise. \\ \\
  \textit{Case 1}:
  $cb(1+\frac{b}{a}) > 0$. In this case if $u > 0$ near the origin that is the equilibrium point, we have $cb(1+\frac{b}{a})u^2 > 0$. If $u < 0$ near the origin we have $cb(1+\frac{b}{a})u^2 > 0$. In both case we have $\dot{u} > 0$ that is the velocity vector of the $u$ direction is along the positive $u$ axis. So the origin that is the equilibrium point is unstable \\ \\
   \textit{Case 2}:
  $cb(1+\frac{b}{a}) < 0$. In this case if $u > 0$ near the origin that is the equilibrium point, we have $cb(1+\frac{b}{a})u^2 < 0$. If $u < 0$ near the origin we have $cb(1+\frac{b}{a})u^2 < 0$. In both case we have $\dot{u} < 0$ that is the velocity vector of the $u$ direction is along the negative $u$ axis. So the origin that is the equilibrium point is unstable \\ \\
  So we see in either of cases \textit{Case 1} and \textit{Case 2}, the origin that is the equilibrium point is unstable.

  \subsubsection{Verification of \textit{Theorem 1}}
  We have the system in the eigenbasis as,
   \begin{equation}
  \begin{split}
  & \dot{u}=\frac{1}{a}\left[c(bu+v)^2+c(bu+v)au\right] \\ &
  \dot{v}=aX_{e}v+a(bu+v)^2-b(bu+v)au-\frac{b}{a}\left[c(bu+v)^2+c(bu+v)au\right]
  \end{split}
  \end{equation}
  Denote,
  \begin{equation}
  \begin{split}
 & f(u,v)=\frac{1}{a}\left[c(bu+v)^2+c(bu+v)au\right] \\ &
 g(u,v)=a(bu+v)^2-b(bu+v)au-\frac{b}{a}\left[c(bu+v)^2+c(bu+v)au\right]
  \end{split}
  \end{equation}
  These are the nonlinear parts of the system. \\
  In the above equations puting $u=-u$ and $v=-v$ we see that 
  \begin{equation}
  \begin{split}
  & f(u,v)=f(-u,-v) \\ &
  g(u,v)=g(-u,-v)
  \end{split}
  \end{equation}
  that is they are even. On the other hand the centre manifold equation is,
  \begin{equation}
  \begin{split}
  h(u)=\frac{b^2c}{aX_{e}}(1+\frac{b}{a})u^2
  \end{split}
  \end{equation}
  which is even and has leading order term as a even power of $u$. This verifies our \textit{Theorem 1} and the conjecture made in \textit{Theorem 1}.
  
  \section{Conclusion}
  Through our work we have discussed the centre manifold theory with application to a general $3-D$ nonlinear system with second order nonlinearities. We did the stability analysis of the system after reducing the state equations on a $2-D$ centre manifold. The reduced phase space was found to possess a stable and unstable fixed point. The system approaches the stable fixed point which is the origin from the unstable fixed point in asymptotic time with a spiral motion towards the origin. In the $x-t$ and $y-t$ plane it was shown to exhibit a oscillatory behaviour with amplitude decaying with time. Then we cast the example of protein molecules which is a $2-D$ nonlinear system. Stability analysis is done around a fixed point which happens to be the origin in our case. The equations obtained point towards a unstable system which very much matches with the predictions of standard theory. 
  With reference\cite{tsuruyama2017kinetic} to the chosen example of the kinetics of protein assembly,  the success of our study lies
in the fact that we have shown mathematically for the oscillation to occur
we require a minimum of three original dimensions with a reduced system of 2-
dimensions on the centre manifold for obtaining the phase plane where the
oscillations take place. If the dimension of the original system is two we can only predict the
stability around the equilibrium point as we have shown in our example.
  
  Along the course of discussion we also state and prove a theorem along with a conjecture made in the theorem. We showed that our theorem and conjecture are true in the light of obtained results for the example that we consider. The significance of the theorem that we prove lies in the fact that we can at once do the stability analysis around the equilibrium point without delving into detailed analysis. Centre manifold analysis requires finding the centre manifold equation and further the reduced system of equations on the centre manifold. But if the nonlinear parts of the decoupled equations are either even or odd we can at once guess the parity of the centre manifold equation and subsequently the stability of the system around the equilibrium point without doing the detailed analysis.\par  
  As a future extension of our work we would like to mention the above analysis could be done along the line shown in our work for a standard $n-D$ system with higher order nonlinearities. Higher dimensional systems often arise in physics for example a Hamiltonian system with two degrees of freedom where the two positions and two momenta give rise to a 4-D system of O.D.E in position and momentum. Higher dimensional systems with high degrees of nonlinearities are also encountered frequently in engineering problems specifically those related to aerodynamical models. So investigating such a system is of natural interest to see if the predictions of the theory also incorporates the higher order nonlinearities of a higher dimensional system and if the theory works well in all different ranges and all different situations.
  \bibliographystyle{unsrtnat}
\bibliography{library}
\end{document}